\documentclass{article}
\usepackage{amsfonts,graphicx,lscape}

\usepackage{amssymb}
\usepackage{eucal}
\usepackage{amsmath}
\usepackage{amsthm}
\usepackage{amscd}
\usepackage{graphics}
\usepackage{graphicx}
\usepackage{amsfonts}
\usepackage{latexsym}
\usepackage[all]{xy}

\numberwithin{equation}{section}
\newtheorem{thm}{\bf Theorem}[section]

\newtheorem{prop}[thm]{\bf Proposition}

\theoremstyle{remark}

\newcommand{\bc}{\begin{center}}
\newcommand{\ec}{\end{center}}
\newcommand{\bec}{\begin{equation}}
\newcommand{\eec}{\end{equation}}
\newcommand{\bea}{\begin{eqnarray}}
\newcommand{\eea}{\end{eqnarray}}
\newcommand{\ba}{\begin{array}}
\newcommand{\ea}{\end{array}}

\def\ds{\displaystyle}
\def\vs{\vspace{4pt}}

\begin{document}

\title{On the symmetries of a Rikitake type system}
\author{Cristian L\u azureanu and Tudor B\^inzar\\
{\small Department of Mathematics, "Politehnica" University of Timi\c soara}\\
{\small Pia\c ta Victoriei nr. 2, 300006 Timi\c soara, Rom\^ania}\\
{\small E-mail: cristian.lazureanu@upt.ro;tudor.binzar@upt.ro}}
\date{}

\maketitle

\begin{abstract}
\noindent A symplectic realization and some
symmetries of a Rikitake type system are presented.
\footnote{COMPTES RENDUS MATHEMATIQUE  Vol. 350,  Issue: 9-10 (2012)  Pages: 529--533;   DOI: 10.1016/j.crma.2012.04.016}
\end{abstract}


\section{Introduction}
The Rikitake two-disk dynamo system is a mechanical model used to study the
reversals of the Earth's magnetic field (\cite{Riki}).

This system has been widely investigated from different points of view. In a particular
case, continuous symmetries were given by W.H.Steeb (\cite{Steeb}).

In this paper, considering another particular case of Rikitake system, namely
\begin{equation}\label{eq1.1}
\left\{\begin{array}{l}
\dot x=yz+\beta y\\
\dot y=xz-\beta x\\
\dot z=-xy
\end{array}\right.,
\end{equation}
where $\beta\in{\bf R}$, some symmetries are given.

A similar study for Maxwell-Bloch equations was presented by P.A.Damianou and P.G.Paschali in
\cite{DamPas}.

Theoretical details about symmetries of differential equations can be found in \cite{BluKum},
\cite{Olver}, \cite{FokFuc}, \cite{Fuc}, \cite{Dam}.

For our purposes, a Hamilton-Poisson realization and a symplectic realization
of system (\ref{eq1.1}) are required.

Symmetry and Hamiltonian systems are related in the topic "order in
chaos". Physical systems often exhibit "order" simultaneous with
symmetries. When the symmetry of a system is broken, Hamiltonian
structures can be useful in detecting chaos \cite{MarWei}.

\section{Bi-Hamiltonian structure and symmetries}
In this section, we consider system (\ref{eq1.1}) with $\beta =0$, i.e.
\begin{equation}\label{eq2.1}
\left\{\begin{array}{l}
\dot x=yz\\
\dot y=xz\\
\dot z=-xy
\end{array}\right.
\end{equation}

A bi-Hamiltonian structure and some symmetries of system (\ref{eq2.1}) are presented.

Let us considering the three-dimensional Lie group of rigid motions of the Minkowski plane,
$$E(1,1)=\{A\in GL(3,{\bf R})|~A=\left[\begin{array}{ccc}
1&0&0\\
v_1&\cosh\theta &\sinh\theta\\
v_2&\sinh\theta &\cosh\theta\end{array}\right],~v_1,v_2,\theta\in{\bf R}\}.$$

The corresponding Lie algebra of $E(1,1)$ is
$$e(1,1)=\{X\in gl(3,{\bf R})|~X=\left[\begin{array}{ccc}
0&0&0\\
a&0&c\\
b&c&0\end{array}\right],~a,b,c\in{\bf R}\}.$$

Note that, as a real vector space, $e(1,1)$ is generated by the base
$B_{e(1,1)}=\{E_1,E_2,E_3\},$ where
$$E_1=\left[\begin{array}{rrr}
0&0&0\\
-\ds{\frac{1}{2}}&0&0\\
0&0&0\end{array}\right]~,~~E_2=\left[\begin{array}{rrr}
0&0&0\\0&0&0\\
-\ds{\frac{1}{2}}&0&0\end{array}\right]~,~~E_3=\left[\begin{array}{rrr}
0&0&0\\0&0&-\ds{\frac{1}{2}}\\
0&-\ds{\frac{1}{2}}&0\end{array}\right].$$

The following bracket relations $[E_1,E_2]=0,[E_1,E_3]=\ds\frac{1}{2}E_2,
[E_2,E_3]=\ds\frac{1}{2}E_1,$ hold.

On the dual space $e(1,1)^*\simeq{\bf R}^3$, the Lie-Poisson structure is given
in coordinates using matrix notation by
$$\pi_1(x,y,z)=\left[\begin{array}{ccc}0&0&\ds\frac{1}{2}y\vs\\
\ds0&0&\ds\frac{1}{2}x\vs\\-\ds\frac{1}{2}y&-\ds\frac{1}{2}x&0\end{array}\right].$$

Following \cite{TudAroNic}, considering the Lie group $O(Q)=\{A\in GL(3,{\bf R})|~
A^tQA=Q\}$ generated by $Q:=diag(2,1,1)\in GL(3,{\bf R})$, the corresponding
Lie algebra is $o(Q)=\{X\in gl(3,{\bf R})|~X^tQ+QX=O_3\}.$ As a real vector space $o(Q)$
is generated by the base $B_{o(Q)}=\{X_1,X_2,X_3\},$ where
$$X_1=\left[\begin{array}{rrr}
0&1&0\\
-\ds{2}&0&0\\0&0&0\end{array}\right]~,~~X_2=\left[\begin{array}{rrr}
0&0&1\\0&0&0\\-\ds{2}&0&0\end{array}\right]~,~~X_3=\left[\begin{array}{rrr}
0&0&0\\0&0&1\\0&-1&0\end{array}\right].$$

The following bracket relations $[X_1,X_2]=-2X_3,[X_1,X_3]=X_2,
[X_2,X_3]=-X_1,$ hold.

On the dual space $o(Q)^*\simeq{\bf R}^3$, the Lie-Poisson structure is given
in coordinates using matrix notation by
$$\pi_2(x,y,z)=\left[\begin{array}{ccc}0&-2z&y\\
2z&0&-x\\-y&x&0\end{array}\right].
$$
\par Taking the constants of motion
$~H_1=\frac{1}{4}x^2-\frac{1}{4}y^2~~\mbox{and}~~H_2=\frac{1}{2}x^2+\frac{1}{2}y^2+z^2~$
of system (\ref{eq2.1}), the following relations
$$\pi_1\cdot\nabla H_2=\pi_2\cdot\nabla H_1=\left(\begin{array}{c}yz\\xz\\-xy\end{array}\right)$$
hold.

Thus, system (\ref{eq2.1}) is a bi-Hamiltonian system. For $\pi_1$ bracket, $H_2$ is the Hamiltonian and
$H_1$ is a Casimir. For $\pi_2$ bracket, $H_2$ is the Hamiltonian and
$H_2$ is a Casimir.

We recall that for a system $\dot x=f(x)$, where $f:M\to TM$, and $M$ is a smooth manifold of finite dimension, a vector
field ${\bf X}$ is called:

$\bullet $ {\em a Lie-point symmetry} if its first prolongation transforms solutions of the system into other solutions;

$\bullet $ {\em a conformal symmetry} if the Lie derivative along ${\bf X}$ satisfies $L_{\bf X}\pi =\lambda\pi $ and
$L_{\bf X}H=\nu H$, for some scalars $\lambda ,\nu $, where the Poisson tensor $\pi $ and the Hamiltonian $H$ give the
Hamilton-Poisson realization of the system;

$\bullet $ {\em a master symmetry} if $[[{\bf X},{\bf X}_f],{\bf X}_f]=0$, but $[{\bf X},{\bf X}_f]\not=0$, where ${\bf X}_f$ is
the vector field defined by the system.

The next result one furnishes a Lie point symmetry of system (\ref{eq2.1}) and a conformal symmetry.
\begin{prop} The vector field
$$~\ds{{\bf X}=-t\frac{\partial }{\partial t}+x\frac{\partial }{\partial x}+
y\frac{\partial }{\partial y}+z\frac{\partial }{\partial z}}~$$
is a Lie point symmetry of system (\ref{eq2.1}). Moreover, ${\bf X}$ is a conformal symmetry.
\end{prop}

{\bf Proof.}  If the vector
$\ds{{\bf v}=\tau (t,x,y,z)\frac{\partial }{\partial t}+A_1(t,x,y,z)\frac{\partial }{\partial x}+
A_2(t,x,y,z)\frac{\partial }{\partial y}+A_3(t,x,y,z)\frac{\partial }{\partial z}}$
is a Lie point symmetry, then its first prolongation
$$pr^{(1)}({\bf v})={\bf v}+(\dot{A}_1-\dot\tau\dot x)\frac{\partial }{\partial \dot x}+
(\dot{A}_2-\dot\tau\dot y)\frac{\partial }{\partial \dot y}+(\dot{A}_3-\dot\tau\dot z)\frac{\partial }{\partial\dot z}$$
applied to our system implies
$$\left\{\begin{array}{l}
\dot{A}_1-yz\dot\tau -zA_2-yA_3=0\\
\dot{A}_2-xz\dot\tau -zA_1-xA_3=0~~.\\
\dot{A}_3+xy\dot\tau +yA_1+xA_2=0\end{array}\right.$$

One solution of above system is the vector ${\bf X}$.

One can easily check that
$$L_{\bf X}\pi_1=-\pi_1~,~~L_{\bf X}\pi_2=-\pi_2~,~~L_{\bf X}H_1=2H_1~,~~
L_{\bf X}H_2=2H_2,$$
whence ${\bf X}$ is a conformal symmetry.

The following result provides a master symmetry of our considered system.

\begin{prop} The vector field
$$~\ds{\overrightarrow{X}=(k_1x+k_2yz)\frac{\partial }{\partial x}+(k_1y+k_2xz)\frac{\partial }{\partial y}+
(k_1z-k_2xy)\frac{\partial }{\partial z}}~,~~k_1\in{\bf R}^*,~k_2\in{\bf R},$$
is a master symmetry of system (\ref{eq2.1}).
\end{prop}

{\bf Proof.} We denote by $\overrightarrow{V}$ the associated
vector field of system (\ref{eq2.1}), that is
$\ds{\overrightarrow{V}=yz\frac{\partial }{\partial
x}+xz\frac{\partial }{\partial y}-xy\frac{\partial }{\partial
z}}.$ It follows that the following relations
$~[\overrightarrow{X},\overrightarrow{V}]=k_1\overrightarrow{V},~\left[[\overrightarrow{X},\overrightarrow{V}],\overrightarrow{V}\right]=0$
hold.

Therefore $\overrightarrow{X}$ is a master symmetry of system (\ref{eq2.1}).

\section{Symplectic realization and symmetries}
Let us consider system (\ref{eq1.1}) in the case $\beta\not=0$.

In this section a symplectic realization of system (\ref{eq1.1}) is given. Using this fact, the symmetries of Newton's equations are studied.

In order to obtain a Hamilton-Poisson realization of system (\ref{eq1.1}),
we again consider the Lie algebra $e(1,1)$ having now the base
$B_{e(1,1)}^{\beta }=\{E_1^\beta ,E_2^\beta ,E_3^\beta\}$, where
$$E_1^\beta =\left[\begin{array}{ccc}
0&0&0\\
-\ds{\frac{1}{2\beta }}&0&0\\
0&0&0\end{array}\right]~,~~E_2^\beta =\left[\begin{array}{ccc}
0&0&0\\0&0&0\\
-\ds{\frac{1}{2\beta }}&0&0\end{array}\right]~,~~E_3^\beta =\left[\begin{array}{ccc}
0&0&0\\0&0&-\ds{\frac{1}{2\beta }}\\
0&-\ds{\frac{1}{2\beta }}&0\end{array}\right].$$
with $[E_1^\beta ,E_2^\beta ]=0,[E_1^\beta ,E_3^\beta ]=
\ds\frac{1}{2\beta }E_2^\beta ,[E_2^\beta ,E_3^\beta ]=
\ds\frac{1}{2\beta }E_1^\beta .$

Following \cite{LibMar}, it is easy to see that the bilinear map
$\Theta :e(1,1)\times e(1,1)\to{\bf R}$ given by the matrix
$(\Theta_{ij})_{1\leq i,j\leq 3},$ $\Theta_{12}=-\Theta_{21}=1$ and
0 otherwise, is a 2-cocycle on $e(1,1)$ and it is not a coboundary since
$\Theta (E_1^\beta ,E_2^\beta )=1\not=0=f([E_1^\beta , E_2^\beta ])$, for
every linear map $f,~f:e(1,1)\to{\bf R}.$

On the dual space $e(1,1)^*\simeq{\bf R}^3$, a modified Lie-Poisson structure
is given in coordinates by
$$\Pi^\beta =\left[\begin{array}{ccc}0&0&\ds\frac{1}{2\beta}y\vs\\
\ds 0&0&\ds\frac{1}{2\beta}x\vs\\
-\ds\frac{1}{2\beta}y&-\ds\frac{1}{2\beta}x&0\end{array}\right]+
\left[\begin{array}{rrr}
0&1&0\\-1&0&0\\
0&0&0\end{array}\right]=
\left[\begin{array}{ccc}0&1&\ds\frac{1}{2\beta}y\vs\\
\ds -1&0&\ds\frac{1}{2\beta}x\vs\\
-\ds\frac{1}{2\beta}y&-\ds\frac{1}{2\beta}x&0\end{array}\right]$$

The Hamiltonian $H_{\beta}$ is given by
$H_{\beta }(x,y,z)=\ds{\frac{\beta }{2}x^2+\frac{\beta }{2}y^2+\beta z^2},$ and
moreover, the function $C_{\beta }$,
$C_{\beta }(x,y,z)=\ds{\frac{1}{4\beta }x^2-\frac{1}{4\beta }y^2+z}$,
is a Casimir of our configuration.

The next theorem states that the system (\ref{eq1.1}) can be regarded as a Hamiltonian mechanical system.

\begin{thm} The Hamilton-Poisson mechanical system $({\bf R}^3,\Pi^\beta ,H_{\beta })$ has a full symplectic realization
$({\bf R}^4,\omega ,H)$, where $\omega =\mbox{d}p_1\wedge\mbox{d}q_1+\mbox{d}p_2\wedge\mbox{d}q_2$ and
$$H=\frac{1}{16\beta }q_1^4+\frac{1}{16\beta }p_1^4-\frac{1}{8\beta }q_1^2p_1^2-\frac{1}{2}q_1^2p_2+\frac{1}{2}p_1^2p_2+
\frac{\beta }{2}q_1^2+\frac{\beta }{2}p_1^2+\beta p_2^2.$$
\end{thm}

{\bf Proof.} The corresponding Hamilton's equations are
\begin{equation}\label{eq3.1}
\left\{\begin{array}{l}
\dot{q}_1=\ds\frac{1}{4\beta }p_1^3-\ds\frac{1}{4\beta }q_1^2p_1+p_1p_2+\beta p_1\vs\\
\dot{q}_2=-\ds\frac{1}{2}q_1^2+\ds\frac{1}{2}p_1^2+2\beta p_2\vs\\
\dot{p}_1=-\ds\frac{1}{4\beta }q_1^3+\ds\frac{1}{4\beta }q_1p_1^2+q_1p_2-\beta q_1\vs\\
\dot{p}_2=0
\end{array}\right.
\end{equation}

We define the application $\varphi :{\bf R}^4\to{\bf R}^3$ by
$$\varphi (q_1,q_2,p_1,p_2)=\left(q_1,p_1,-\ds\frac{1}{4\beta }q_1^2+\ds\frac{1}{4\beta }p_1^2+p_2\right)=(x,y,z).$$

It follows that $\varphi $ is a surjective submersion, the equations (\ref{eq3.1}) are mapped onto the equations (\ref{eq1.1}), the
canonical structure $\{.,.\}_{\omega }$ is mapped onto the Poisson structure $\Pi^\beta $, as required.

We also remark that $H_{\beta }=H$ and $C_{\beta }=p_2.$\\

From Hamilton's equations (\ref{eq3.1}) we obtain by differentiation, Newton's equations:
$$\ddot{q}_2\dot{q}_2+2\beta^2\ddot{q}_2+4\beta^2q_1\dot{q}_1=0$$
$$2\beta \ddot{q}_1\dot{q}_2+4\beta^3\ddot{q}_1-2\beta\dot{q}_1\ddot{q}_2-\frac{1}{2\beta }q_1\dot{q}_2^3-\beta q_1\dot{q}_2^2+
2\beta^3 q_1\dot{q}_2+4\beta^5q_1=0$$
These are also Lagrange's equations generated by the Lagrangian
$$L=\frac{1}{4\beta }\dot{q}_2^2-\frac{\beta }{2}q_1^2+\frac{1}{4\beta }q_1^2\dot{q}_2+\frac{\beta\dot{q}_1^2}{\dot{q}_2+2\beta^2}.$$

A vector field $$\overrightarrow{v}=\xi (q_1,q_2,t)\frac{\partial }{\partial t}+\eta_1(q_1,q_2,t)\frac{\partial }{\partial q_1}+
\eta_2(q_1,q_2,t)\frac{\partial }{\partial q_2}$$
is a Lie Point symmetry for Newton's equations if the action of its second prolongation on Newton's equations vanishes. Thus, the
following conditions are obtained:\\
$(\dot{q}_2+2\beta^2)\ddot{\eta}_2-\dot{q}_2(\dot{q}_2+2\beta^2)\ddot{\xi }+4\beta^2q_1\dot{\eta }_1+\ddot{q}_2\dot{\eta }_2-
\dot{\xi }(4\beta^2q_1\dot{q}_1+3\ddot{q}_2\dot{q}_2+4\beta^2\ddot{q}_2)+4\beta^2\dot{q}_1\eta_1=0$\\
$2\beta (\dot{q}_2+2\beta^2)\ddot{\eta}_1-2\beta\dot{q}_1\ddot{\eta }_2-4\beta^3\dot{q}_1\ddot{\xi }-2\beta\ddot{q}_2\dot{\eta }_1+
(2\beta\ddot{q}_1-\frac{3}{2\beta }q_1\dot{q}_2^2-2\beta q_1\dot{q}_2+2\beta^3q_1)\dot{\eta }_2+
(6\beta\dot{q}_1\ddot{q}_2-6\beta\ddot{q}_1\dot{q}_2-8\beta^3\ddot{q}_1+\frac{3}{2\beta }q_1\dot{q}_2^3+
2\beta q_1\dot{q}_2^2-2\beta^3q_1\dot{q}_2)\dot{\xi }+(-\frac{1}{2\beta }\dot{q}_2^3-\beta\dot{q}_2^2+2\beta^3\dot{q}_2+4\beta^5)\eta_1=0.$

The resulting equations obtained by expanding $\dot{\xi },\ddot{\xi },\dot{\eta }_1,\ddot{\eta }_1,\dot{\eta }_2,\ddot{\eta }_2$ must be
satisfied identically in $t,$ $q_1,$ $q_2,$ $\dot{q}_1,$ $\dot{q}_2$, that are all independent. Doing standard manipulation, we get the overall result:
$$\left\{\begin{array}{l}\xi =c_1\\\eta_1=0\\\eta_2=c_2\end{array}\right.$$
where $c_1,c_2$ are real constants.

For $c_2=0$, but $c_1\not=0$, we have $\overrightarrow{v_1}=c_1\ds\frac{\partial }{\partial t}$ that represents the time translation symmetry
which generates the conservation of energy $H$.

For $c_1=0$, but $c_2\not=0$, we have $\overrightarrow{v_2}=c_2\ds\frac{\partial }{\partial q_2}$ that represents a translation in the
cyclic $q_2$ direction which is related to the conservation of $p_2$.

Moreover, using the Lagrangian $L$ and Noether's theory we deduce that both $\overrightarrow{v_1}$ and $\overrightarrow{v_2}$ are variational symmetries
since they satisfy the condition $~pr^{(1)}\overrightarrow{v}(L)+L\mbox{div}~(\xi )=0.$

\section*{Acknowledgements}
We would like to thank the referees very much for
their valuable comments and suggestions.

\end{document}